\documentclass[aos,seceqn,dvips]{arximspdf}
\usepackage{mathbh}
\usepackage{graphicx}

\doi{10.1214/10-AOS828}
\volume{39}
\issue{2}
\pubyear{2011}
\firstpage{990}
\lastpage{1011}

\makeatletter
\newcommand{\rset}{\mathbb R}
\renewcommand{\epsilon}{\varepsilon}
\newcommand{\un}{\mathbh{1}}
\newcommand{\eqdef}{\ensuremath{\stackrel{\mathrm{def}}{=}}}

\newcommand{\Xset}{\mathsf{X}} 
\newcommand{\Xsigma}{\mathcal{X}} 
\newcommand{\Tsigma}{\mathcal{B}(\Theta)} 
\newcommand{\F}{\mathcal{F}} 

\newcommand{\PP}{\mathbb{P}} 
\newcommand{\PE}{\mathbb{E}} 

\newcommand{\Cset}{\mathcal{C}} 

\newcommand{\p}{\mathbb{P}}
\newcommand{\E}{\mathbb{E}}
\newcommand{\Real}{\mathbb{R}}

\newtheorem{theo}{Theorem}[section]
\newtheorem{coro}{Corollary}[section]
\newtheorem{proposition}{Proposition}[section]
\newproclaim{algo}{Algorithm}[section]
\newtheorem{lemma}{Lemma}[section]

\newproclaim{definition}{Definition}[section]
\newproclaim{example}{Example}[section]
\newproclaim{remark}{Remark}[section]

\makeatother

\begin{document}
\begin{frontmatter}

\title{Kernel estimators of asymptotic variance\\ for adaptive Markov
chain Monte Carlo\thanksref{T1}}
\runtitle{Asymptotic variance estimation for adaptive MCMC}
\thankstext{T1}{Supported in part by NSF Grant DMS-09-06631.}

\begin{aug}
\author[A]{\fnms{Yves F.} \snm{Atchad\'e}\corref{}\ead[label=e1]{yvesa@umich.edu}}
\runauthor{Y. F. ATCHAD\'E}
\affiliation{University of Michigan}
\address[A]{Department of Statistics\\
University of Michigan\\
1085 S. University Avenue\\
Ann Arbor, Michigan 48109 \\
USA\\
\printead{e1}} 
\end{aug}

\received{\smonth{11} \syear{2009}}
\revised{\smonth{4} \syear{2010}}

%
\begin{abstract}
We study the asymptotic behavior of kernel estimators of asymptotic
variances (or long-run variances) for a class of adaptive Markov
chains.
The convergence is studied both in $L^p$ and almost surely. The
results also apply to Markov chains and improve on the existing
literature by imposing weaker conditions. We illustrate the results
with applications to the $\operatorname{GARCH}(1,1)$ Markov model and to an adaptive
MCMC algorithm for Bayesian logistic regression.
\end{abstract}

%

\end{frontmatter}

\section{Introduction}
Adaptive Markov chain Monte Carlo (adaptive MCMC) provides a flexible
framework for optimizing MCMC samplers on the fly (see, e.g.,
\cite{andrieuetthoms08,atchadeetal09,robertsetrosenthal06} and the
reference therein). If $\pi$ is the probability measure of interest,
then these adaptive MCMC samplers generate random processes
$\{X_n,n\geq0\}$ that typically are not Markov, but they
nevertheless satisfy a~law of large numbers and the empirical average
$n^{-1}\sum_{k=1}^n h(X_k)$ provides a consistent estimate of the
integral $\pi(h)\eqdef\PE(h(X))$, $X\sim\pi$. A measure of
uncertainty in approximating $\pi(h)$ by the random variable
$n^{-1}\sum_{k=1}^n h(X_k)$ is given by the variance
$\operatorname{Var}(n^{-1/2}\sum_{k=1}^nh(X_k))$. In particular, the
asymptotic variance $\sigma^2(h)\eqdef
\lim_{n\to\infty}\operatorname{Var}(n^{-1/2}\sum_{k=1}^nh(X_k))$
(also known as the \textit{long-run variance}) plays a fundamental role
in assessing the performances of Monte Carlo simulations. But the
problem of estimating asymptotic variances for adaptive MCMC samplers
has not been addressed in the literature.

We study kernel estimators of asymptotic variances for a general class
of adaptive Markov chains. These adaptive Markov chains (the precise
definition is given in Section \ref{sec:adaptchain} below), which
include Markov chains, constitute a~theoretical framework for analyzing
adaptive MCMC algorithms. More precisely, if $\{X_n,n\geq0\}$ is an
adaptive Markov chain and $h\dvtx\Xset\to\Real$ a function of interest,
then we consider estimators of the form
\[
\Gamma^2_n(h)=\sum_{k=-n}^nw(k b)\gamma_n(k),
\]
where $\gamma_n(k)=\gamma_n(k;h)$ is the $k$th order sample
autocovariance of $\{h(X_n),n\geq0\}$, $w\dvtx\Real\to\Real$ is a
kernel with support $[-1,1]$ and $b=b_n$ is the bandwidth. These are
well-known methods pioneered by M. S. Bartlett, M. Rosenblatt, E.
Parzen and others (see, e.g., \cite{priestley81} for more details).
But, with a few notable exceptions in the econometrics literature (see
references below), these estimators have mostly been studied with the
assumption of stationarity. Thus, more broadly, this paper contributes
to the literature on the behavior of kernel estimators of asymptotic
variances for ergodic nonstationary processes.

It turns out that, in general, the asymptotic\vspace*{1pt} variance
$\sigma^2(h)$ does not characterize the limiting distribution of
$n^{-1/2}\sum_{k=1}^n (h(X_k)-\pi(h))$ as, for example, with ergodic
Markov chains. For adaptive Markov chains, we show that
$n^{-1/2}\sum_{k=1}^n (h(X_k)-\pi(h))$ converges weakly to a mixture of
normal distributions of the form $\sqrt{\Gamma^2(h)}Z$ for some mixing
random variable $\Gamma^2(h)$, where $Z$ is a standard normal random
variable independent of $\Gamma^2(h)$.
Under a geometric drift stability condition on the adaptive Markov
chain and some verifiable conditions on the kernel $w$ and the
bandwidth $b_n$, we prove that the kernel estimator $\Gamma^2_n(h)$
converges to $\Gamma^2(h)$ in $L^p$-norm, $p>1$, and almost surely.
For Markov chains, $\Gamma^2(h)$ coincides with $\sigma^2(h)$, the
asymptotic variance of $h$. Another important special case where we
have $\Gamma^2(h)=\sigma^2(h)$ is the one where the adaptation
parameter converges to a deterministic limit as, for instance, with the
adaptive Metropolis algorithm of \cite{haarioetal00}. The general case
where $\Gamma^2(h)$ is random poses some new difficulties to Monte
Carlo error assessment in adaptive MCMC that we discuss in Section
\ref{finalrem}.

We derive the rate of convergence for $\Gamma_n^2(h),$
which suggests\vspace*{1pt} selecting the bandwidth to be $b_n\propto
n^{-(2/3)(1-0.5\vee(1/p))}$. When $p=2$ is admissible, we obtain the
bandwidth $b_n\propto n^{-1/3}$, as in \cite{jonesetal09}.

The problem of estimating asymptotic variances is well known in MCMC
and Monte Carlo simulation in general. Besides the estimator described
above, several other methods have been proposed, including batch means,
overlapping batch means and regenerative simulation
(\cite{bratleyetal87,damerdji95,jonesetal09,myklandetal95}). For the
asymptotics of kernel estimators, the important work of
\cite{jonesetal09} proves the $L^2$-consistency and strong consistency
of kernel estimators for Markov chains under the assumption of
geometric ergodicity and $\E(|h(X)|^{4+\epsilon})<\infty$, $X\sim
\pi$,
for some $\epsilon>0$. We weaken these moment conditions to
$\E(|h(X)|^{2+\epsilon})<\infty$.

Estimating asymptotic variances is also a well-known problem in
econometrics and time series modeling. For example, if $\hat\beta_n$ is
the ordinary least-squares estimator of $\beta$ in the simple linear
model $y_i=\alpha+ \beta x_i +u_i,i=1,\ldots,n,$ where
$\{u_k,k\geq1\}$ is a dependent noise process, then, under some mild
conditions on the sequence $\{x_i\}$ and on the noise process,
$\sqrt{n}(\hat\beta_n-\beta)$ converges weakly to a normal
distribution  $\mathcal{N}(0,\sigma^2/c^2),$ where

{\fontsize{10.3}{11.3}\selectfont{
\[\sigma^2=\lim_{n\to\infty}\!\operatorname{Var}\Biggl(\!n^{-1/2}\sum_{k=1}^n
u_k\!\Biggr),\qquad
c^2=\lim_{n\to\infty}n^{-1}\!\sum_{k=1}^n(x_i-\bar x_n)^2,\qquad
\bar x_n=n^{-1}\!\sum_{k=1}^n x_k.
\]
}}
Therefore, a valid inference on $\beta$ requires the estimation of the
asymptotic variance $\sigma^2$. The multivariate version of this
problem involves estimating the so-called \textit{heteroskedasticity
and autocorrelation} (HAC) matrices. Several authors have studied the
kernel estimation of HAC matrices and attention has been paid to
nonstationarity under various mixing assumptions or mixingale-type
assumptions (\cite{andrew91,dejong00,dejongetdavidson00,hansen92}).
But these results require mixing conditions that do not hold in the
present setup.

On a more technical note, the proof of our main results (Theorems
\ref{thm1}--\ref{thm2}) is based on a martingale approximation approach
adapted from \cite{wuetshao07}. The crux of the argument consists in
approximating the periodogram of the adaptive Markov chain by a
quadratic form of a martingale difference process which is then treated
as a martingale array. As part of the proof, we develop a strong law of
large numbers for martingale arrays which may also be of some
independent interest. The approach taken here thus differs from the
almost sure strong approximation approach taken in
\cite{damerdji95,jonesetal09}.

The paper is organized as follows. In Section \ref{sec:adaptchain}, we
define the class of adaptive Markov chains that will be studied. In
Section \ref{sec:clt}, we give a general central limit theorem for
adaptive Markov chains that sets the stage to better understand the
limiting behavior of the kernel estimator $\Gamma^2_n(h)$.
In Section~\ref{sec:kernelmethod}, we state the assumptions and the main results
of the paper. We also discuss some practical implications of these
theoretical results. The proofs are postponed to Section
\ref{sec:proofs} and to the supplementary paper~\cite{atchadesuppl}.
Section \ref{sec:examples} presents applications to generalized
autoregressive conditional heteroscedastic ($\operatorname{GARCH}$) processes and to a
Bayesian analysis of logistic regression.

We end this introduction with some general notation that will be used
throughout the paper. For a Markov kernel $Q$ on a measurable space
$(\mathcal{Y},\mathcal{A}),$ say, we denote by $Q^n$, $n\geq0$, its
$n$th iterate. Any such Markov kernel $Q$ acts both on bounded
measurable functions $f$ and on $\sigma$-finite measures $\mu$, as in
$Qf(\cdot) \eqdef\int Q(\cdot,dy) f(y)$ and $\mu Q(\cdot) \eqdef
\int
\mu(dx) Q(x, \cdot)$. If $W\dvtx \mathcal{Y}\to[1, +\infty)$ is a
function, then the $W$-norm of a function $f\dvtx \mathcal{Y}\to\Real$ is
defined as $|f|_W \eqdef\sup_{\mathcal{Y}} |f| /W$. The set of
measurable functions $f\dvtx \mathcal{Y}\to\Real$ with finite $W$-norm is
denoted by $\mathcal{L}_W$. Similarly, if $\mu$ is a signed measure on
$(\mathcal{Y},\mathcal{A})$, then the $W$-norm of $\mu$ is defined as
$\| \mu\|_{W} \eqdef\sup_{\{g, |g|_W \leq1 \}} |\mu(g)|$, where
$\mu(g)\eqdef\int g(y)\mu(dy)$. If $\nu$ \vspace*{1pt}is a $\sigma$-finite measure
on $(\mathcal{Y},\mathcal{A})$ and $q\geq1$, we denote by $L^q(\nu)$
the space of all measurable functions
$f\dvtx(\mathcal{Y},\mathcal{A})\to\rset$ such that $\nu(|f|^q)<\infty$.
Finally, for $a,b\in\rset$, we define $a\wedge b=\min(a,b)$ and
$a\vee
b=\max(a,b)$.

\section{Adaptive Markov chains}\label{sec:adaptchain}
Let $(\Xset,\Xsigma)$ be a measure state space measure space endowed
with a countably generated $\sigma$-field $\Xsigma$. Let
$(\Theta,\Tsigma)$ be a measure space. In practice, we will take
$\Theta$ to be a compact subspace of $\rset^q$, the $q$-dimensional
Euclidean space. Let $\{P_\theta, \theta\in\Theta\}$ be a family of
Markov transition kernels on $(\Xset,\Xsigma)$ such that for any $(x,A)
\in\Xset\times\Xsigma$, $\theta\mapsto P_\theta(x,A)$ is
measurable. Let $\pi$ be a probability measure on $(\Xset,\Xsigma)$.
We assume that for each $\theta\in\Theta$, $P_\theta$ admits $\pi$ as
its invariant distribution.

The stochastic processes of interest in this work are defined as
follows. Let $\Omega=(\Xset\times\Theta)^\infty$ be the
product space equipped with its product $\sigma$-algebra $\F$ and let
$\bar\mu$ be a probability measure on
$(\Xset\times\Theta,\Xsigma\times\Tsigma)$. Let $\PP_{\bar\mu}$ be the
probability measure on $(\Omega,\F)$ with associated expectation
operator $\PE_{\bar\mu}$, associated process $\{ (X_n,\theta_n), n
\geq0 \}$ and associated natural filtration $\{ \F_n, n \geq0\},$
with the following properties: $(X_0,\theta_0)\sim\bar\mu$ and, for
each $n\geq0$ and any nonnegative measurable function
$f\dvtx\Xset\to\rset$,
%
\begin{equation}\label{eq:condition-markov}
\PE_{\bar\mu}(f(X_{n+1})\vert\F_n)= P_{\theta_n} f(X_n)=
\int P_{\theta_n}(X_n, dy) f(y),\qquad\PP_{\bar\mu}\mbox{-a.s.}
\end{equation}
We call the $\Xset$-marginal process $\{X_n,n\geq0\}$ an
\textit{adaptive Markov chain}. In this definition, we have left the
adaptation dynamics (i.e., the conditional distribution of
$\theta_{n+1}$ given $\F_n$ and $X_{n+1}$) unspecified. This can be
done in many different ways (see, e.g., \cite{robertsetrosenthal06}).
But it is well known, as we will see later, that the adaptation
dynamics needs to be \textit{diminishing} in order for the adaptive
Markov chain to maintain $\pi$ as its limiting distribution.

The simplest example of an adaptive Markov chain is the case where
$\theta_{n}\equiv\bar\theta\in\Theta$ for all $n\geq0$. Then
$\{X_n,n\geq0\}$ is a Markov chain with transition kernel
$P_{\bar\theta}$. In other words, our analysis also applies to Markov
chains and, in particular, to Markov chain Monte Carlo.

\begin{example}\label{exRWM}
To illustrate the definitions and, later, the results, we present a
version of the adaptive Metropolis algorithm of \cite{haarioetal00}.
We take $\Xset=\Real^d$ equipped with its Euclidean norm and inner
product, denoted by $|\cdot|$ and $\langle\cdot,\cdot\rangle$, respectively.
Let $\pi$ be a positive, possibly unnormalized, density (with respect
to the Lebesgue measure). We construct the parameter space $\Theta$ as
follows. We equip the set $\mathcal{M}_+$ of all $d$-dimensional
symmetric positive semidefinite matrices with the Frobenius norm
$|A|\eqdef\sqrt{\operatorname{Tr}(A^T A)}$ and inner product
$\langle A,B\rangle=\operatorname{Tr}(A^T B)$. For $r>0$, let $\Theta_+(r)$ be the
compact subset of elements $A\in\mathcal{M}_+$ such that $|A| \leq r$.
Let $\Theta_\mu(r)$ be the ball centered at $0$ and with radius $r$ in
$\Real^d$. We then define $\Theta\eqdef\Theta_\mu(r_1) \times
\Theta_+(r_2)$ for some constants $r_1,r_2>0$.

We introduce the functions $\Pi_\mu\dvtx\Real^d\to\Theta_\mu(r_1)$ and
$\Pi_+\dvtx\mathcal{M}_+\to\Theta_+(r_2),$ defined as follows. For
$v\in\Theta_\mu(r_1)$, $\Pi_\mu(v)=v$ and for $v\notin\Theta_\mu(r_1)$,
$\Pi_\mu(v)=\frac{M}{|v|}v$. Similarly, for $\Sigma\in\Theta_+(r_2)$,
$\Pi_+(\Sigma)=\Sigma$ and for $\Sigma\notin\Theta_+(r_2)$,
$\Pi_+(\Sigma)=\frac{M}{|\Sigma|}\Sigma$.

For $\theta=(\mu,\Sigma)\in\Theta$, let $P_\theta$ be the transition
kernel of the random walk Metropolis (RWM) algorithm with proposal
kernel $\mathcal{N}(x,\frac{2.38^2}{d}\Sigma+\epsilon I_d)$
and target distribution $\pi$. The adaptive Metropolis algorithm works
as follows.

\begin{algo}\label{arwm1}
\textit{Initialization}: Choose $X_0\in\Real^d$,
$(\mu_0,\Sigma_0)\in\Theta$. Let $\{\gamma_n\}$ be a sequence of
positive numbers (we use $\gamma_n= n^{-0.7}$ in the simulations).

\textit{Iteration}: Given $(X_n,\mu_n,\Sigma_n)$:
\begin{longlist}
\item[(1)] generate $Y_{n+1}\sim
\mathcal{N}(X_n,\frac{2.38^2}{d}\Sigma_n+\epsilon I_d)$;
with probability $\alpha_{n+1}=\alpha(X_n,\break Y_{n+1}),$ set
$X_{n+1}=Y_{n+1}$ and with probability $1-\alpha_{n+1}$, set
$X_{n+1}=X_n$;
\item[(2)] set
%
\begin{eqnarray}\label{ex2:defiMu}
\mu_{n+1} &=& \Pi_\mu\bigl(\mu_n+(n+1)^{-1}(X_{n+1}-\mu_n)\bigr),\\ \label{ex2:defiSigma}
\Sigma_{n+1}
&=&\Pi_+\bigl(\Sigma_n+(n+1)^{-1}\bigl((X_{n+1}-\mu_n)(X_{n+1}-\mu_n)^T-\Sigma_n\bigr)\bigr).
\end{eqnarray}
\end{longlist}

Thus, given $\F_n=\sigma\{X_k,\mu_k,\Sigma_k,k\leq n\}$,
$X_{n+1}\sim
P_{\theta_n}(X_n,\cdot)$, where $P_{\theta_n}$ is the Markov kernel of
the random walk Metropolis with target $\pi$ and proposal
$\mathcal{N}(x,\frac{2.38^2}{d}\Sigma_n+\epsilon I_d)$. So,
this algorithm generates a random process $\{(X_n,\theta_n),n\geq
0\}$ that is an adaptive Markov chain, as defined above. Here, the
adaptation dynamics is given by (\ref{ex2:defiMu}) and
(\ref{ex2:defiSigma}).
\end{algo}
\end{example}

Throughout the paper, we fix the initial measure of the process to some
arbitrary measure $\bar\mu$ and simply write $\E$ and $\PP$ for
$\E_{\bar\mu}$ and $\PP_{\bar\mu}$, respectively. We impose the
following geometric ergodicity assumption.

\begin{longlist}
\item[A1:] For each $\theta\in\Theta$, $P_\theta$ is phi-irreducible
and aperiodic with invariant\break distribution $\pi$. There exists a
measurable function $V\dvtx \Xset\to[1,\infty)$ with\break $\int
V(x)\bar\mu(dx, d\theta)<\infty$ such that for any $\beta\in(0,1]$,
there exist $\rho\in(0,1)$, $C\in(0,\infty)$ such that for any
$(x,\theta) \in\Xset\times\Theta$,
%
\begin{equation}\label{rateconv}
\| P_\theta^n(x,\cdot) - \pi(\cdot) \|_{V^\beta} \leq C \rho^n
V^{\beta}(x),\qquad n\geq0.
\end{equation}
Furthermore, there exist constants $b\in(0,\infty), \lambda\in(0,1)$
such that for any $(x,\theta) \in\Xset\times\Theta$,
%
\begin{equation}\label{drift}
P_\theta V(x)\leq\lambda V(x)+b.
\end{equation}
\end{longlist}

Condition (\ref{rateconv}) is a standard geometric
ergodicity assumption. We impose (\ref{drift}) in order to control the
moments of the adaptive process. Condition (\ref{drift}) is probably
redundant since geometric ergodicity intuitively implies a drift
behavior of the form (\ref{drift}). But this is rarely an issue because
both (\ref{rateconv}) and (\ref{drift}) are implied by the following
minorization and drift conditions.

\begin{longlist}
\item[DR:]Uniformly for
$\theta\in\Theta$, there exist $\Cset\in\Xsigma$, $\nu$ a probability
measure on $(\Xset,\Xsigma)$, $b,\epsilon>0$ and $\lambda\in(0,1)$
such that $\nu(\Cset)>0$, $P_\theta(x,\cdot)\geq
\epsilon\nu(\cdot)\un_\Cset(x)$ and
%
\begin{equation}\label{drift2}
P_\theta V\leq\lambda V+b\un_\Cset.
\end{equation}
\end{longlist}

This assertion follows from Theorem 1.1 of \cite{baxendale05}. DR is
known to hold for many Markov kernels used in MCMC simulation (see,
e.g., \cite{jonesetal09} for some references). Either drift condition
(\ref{drift}) or (\ref{drift2}) implies that $\pi(V)<\infty$
(\cite{meynettweedie93}, Theorem 14.3.7). Therefore, under A1, if
$f\in\mathcal{L}_{V^\beta}$ for some $\beta\in[0,1]$, then $f\in
L^{1/\beta}(\pi)$. Finally, we note that under A1, a law of
large numbers can be established for the adaptive chain (see, e.g.,
\cite{atchadeetfort08}). A short proof is provided here for
completeness.

To state the law of large numbers, we need the following pseudo-metric
on $\Theta$. For $\beta\in[0,1]$, $\theta,\theta'\in\Theta$, set
\[
D_\beta(\theta,\theta') \eqdef\sup_{|f|_{V^\beta}\leq1}\sup_{x
\in
\Xset} \frac{|P_{\theta}f(x) - P_{\theta'}f(x)
|}{V^\beta(x)} .
\]

\begin{proposition}\label{wlln}
Assume \textup{A1}. Let $\beta\in[0,1)$ and
$\{h_\theta\in\mathcal{L}_{V^\beta},\theta\in\Theta\}$ be a family of
functions such that $\pi(h_\theta)=0$, $(x,\theta)\to h_\theta(x)$ is
measurable and $\sup_{\theta\in\Theta}|h_\theta|_{V^\beta}<\infty$.
Suppose also that
%
\begin{equation}\label{diminishwlln}
\sum_{k\geq1}k^{-1}\bigl(D_\beta(\theta_k,\theta_{k-1})
+|h_{\theta_k}-h_{\theta_{k-1}}|_{V^\beta}\bigr)V^{\beta}(X_k)<\infty
,\qquad \PP\mbox{-a.s.}
\end{equation}
Then $n^{-1}\sum_{k=1}^nh_{\theta_{k-1}}(X_k)$ converges almost surely
($\PP$) to zero.
\end{proposition}
\begin{pf}
See Section \ref{proof:propslln}.
\end{pf}

\section{A central limit theorem}\label{sec:clt}
Central limit theorems are useful in assessing Monte Carlo errors.
Several papers have studied central limit theorems for adaptive MCMC
(\cite{andrieuetal06,atchadeetfort08,saksmanvihola09}). The next
proposition is adapted from \cite{atchadeetfort09}. For $h\in\mathcal{L}_V$, we
introduce the resolvent functions
\[
g_\theta(x)\eqdef\sum_{j\geq0}\bar P_{\theta}^j h(x),
\]
where $\bar P_\theta\eqdef P_\theta-\pi$. The dependence of
$g_\theta$
on $h$ is omitted for notational convenience. We also define
$G_\theta(x,y)=g_\theta(y)-P_\theta g_\theta(x)$, where $P_\theta
g_\theta(x)\eqdef\int P_\theta(x,dz)g_\theta(z)$. Whenever
$g_\theta$
is well defined, it satisfies the so-called \textit{Poisson equation}
%
\begin{equation}\label{poisson}
h(x)=g_\theta(x)-\bar P_\theta
g_\theta(x).
\end{equation}

\begin{proposition}\label{propclt}
Assume \textup{A1}. Let $\beta\in[0,1/2)$ and $h\in\mathcal{L}_{V^\beta}$ be such
that $\pi(h)=0$. Suppose that there exists a nonnegative random
variable $\Gamma^2(h),$ finite $\PP$-a.s., such that
%
\begin{equation}\label{asympvar}\lim_{n\to\infty}\frac{1}{n}\sum
_{k=1}^n G_{\theta_{k-1}}^2(X_{k-1},X_k)=\Gamma^2(h)\qquad\mbox{in } \PP\mbox{-probability}.
\end{equation}
Suppose also that
%
\begin{equation}\label{diminishclt}
\sum_{k\geq1}k^{-1/2}D_{\beta}(\theta_k,\theta_{k-1})
V^{\beta}(X_k)<\infty,\qquad\PP\mbox{-a.s.}
\end{equation}
Then
$n^{-1/2}\sum_{k=1}^nh(X_k)$ converges weakly to a random variable
$\sqrt{\Gamma^2(h)}Z$, where $Z\sim\mathcal{N}(0,1)$ is a standard
normal random variable independent of~$\Gamma^2(h)$.
\end{proposition}
\begin{pf}
See Section \ref{proof:propclt}.
\end{pf}

Condition (\ref{diminishclt}), which strengthens (\ref{diminishwlln}),
is a \textit{diminishing adaptation condition} and is not hard to check
in general. It follows from the following assumption which is much
easier to check in practice.

\begin{longlist}
\item[A2:]
There exist $\eta\in
[0,1/2)$ and a nonincreasing
sequence of positive numbers $\{\gamma_n,n\geq1\}$,
$\gamma_n=O(n^{-\alpha})$, $\alpha>1/2$, such that
for any $\beta\in[0,1]$, there exists a finite constant $C$ such that
%
\begin{equation}\label{diminish}
D_\beta(\theta_{n-1},\theta_n)\leq
C\gamma_{n}V^\eta(X_n),\qquad\PP\mbox{-a.s.}
\end{equation}
\end{longlist}
\cite{andrieuetal06} establishes A2 for the random walk
Metropolis and the independence sampler. A similar result is obtained
for the Metropolis adjusted Langevin algorithm in \cite{atchade05}. The
constant $\eta$ in A2 reflects the additional fluctuations due
to the adaptation. For example, for a Metropolis algorithm with
adaptation driven by a stochastic approximation of the form
$\theta_{n+1}=\theta_n+\gamma_n H(\theta_n,X_{n+1})$, $\eta$ is any
nonnegative number such that
$\sup_{\theta\in\Theta}|H(\theta,\cdot)|_{V^\eta}<\infty$.

\begin{proposition}
Under \textup{A1}--\textup{A2}, (\ref{diminishclt}) holds.
\end{proposition}
\begin{pf}
Under A2, the left-hand side of (\ref{diminishclt}) is bounded
almost surely by $C\sum_{k\geq1}k^{-1/2}\gamma_kV^{\eta+\beta}(X_k)$,
the expectation of which is bounded by the term $C\sum_{k\geq
1}k^{-1/2}\gamma_k$ according to Lemma \ref{lemuseful}(a), assuming
A1. Since $\alpha>1/2$, we conclude that (\ref{diminishclt})
holds.
\end{pf}

Equation (\ref{asympvar}) is also a natural assumption. Indeed, in most adaptive
MCMC algorithms, we seek to find the ``best'' Markov kernel from the
family $\{P_\theta,\theta\in\Theta\}$ to sample from $\pi$. Thus, it
is often the case that $\theta_n$ converges to some limit
$\theta_\star$, say (see, e.g.,
\cite{andrieuetal06,andrieuetthoms08,atchadeetfort09,atchadeetrosenthal03}).
In these cases, (\ref{asympvar}) actually holds.

\begin{proposition}\label{propconvtheta}
Assume \textup{A1}--\textup{A2}. Let $\beta\in[0,(1-\eta)/2)$, where
$\eta$
is as in \textup{A2}, and let $h\in\mathcal{L}_{V^\beta}$ be such that $\pi(h)=0$.
Suppose that there exists a~$\Theta$-valued random variable
$\theta_\star$ such that
$D_{\beta}(\theta_n,\theta_\star)+D_{2\beta}(\theta_n,\theta
_\star)$
converges in probability to zero. Then (\ref{asympvar}) holds.
Furthermore,
\[
\Gamma^2(h)=\int_{\Xset\times\Xset}\pi(dx)P_{\theta_\star
}(x,dy)G^2_{\theta_\star}(x,y).
\]
\end{proposition}
\begin{pf}
See Section \ref{proof:convtheta}.
\end{pf}

\begin{definition}\label{def}
We call the random variable $\Gamma^2(h)$ the \textit{asymptotic
average squared variation} of $h$ and $\sigma^2(h)\eqdef
\E(\Gamma^2(h))$ the \textit{asymptotic variance} of $h$.
\end{definition}

This definition is justified by the following result.
\begin{proposition}\label{propvar}Assume \textup{A1}--\textup{A2}. Let
$\beta\in[0,1/2)$ and $h\in\mathcal{L}_{V^\beta}$ be such
that $\pi(h)=0$. Assume that (\ref{asympvar}) holds. Then
\[
\lim_{n\to\infty}\operatorname{Var}\Biggl(n^{-1/2}\sum_{k=1}^nh(X_k)\Biggr)=\sigma^2(h).
\]
\end{proposition}
\begin{pf}
See Section \ref{proof:propvar}.
\end{pf}

\section{Asymptotic variance estimation}\label{sec:kernelmethod}
Denote by $\pi_n(h)=n^{-1}\sum_{k=1}^nh(X_k)$ the sample mean of
$h(X_k)$ and denote by $\gamma_n(k)$ the sample autocovariance:
$\gamma_n(k)=0$ for $|k|\geq n$, $\gamma_n(-k)=\gamma_n(k)$ and for
$0\leq k<n$,
\[
\gamma_n(k)=\frac{1}{n}\sum_{j=1}^{n-k}\bigl(h(X_j)-\pi
_n(h)\bigr)\bigl(h(X_{j+k})-\pi_n(h)\bigr).
\]

Let $w\dvtx\Real\to\Real$ be a function with support $[-1,1]$ [$w(x)=0$
for $|x|\geq1$]. We assume that $w$ satisfies the following.
\begin{longlist}
\item[A3.] The function $w$ is even [$w(-x)=w(x)$] and $w(0)=1$.
Moreover, the restriction $w\dvtx [0,1]\to\Real$ is twice continuously
differentiable.
\end{longlist}

Typical examples of kernels that satisfy A3 include, among
others, the family of kernels
%
\begin{equation}\label{bartlettsK}
w(x)=\cases{1-|x|^q,&\quad if $|x|\leq1$,\cr
 0, &\quad if $|x|>1$,
}
\end{equation}
for $q\geq1$. The case $q=1$ corresponds to the Bartlett kernel.
A3 is also satisfied by the Parzen kernel
%
\begin{equation}\label{parzensK}
w(x)=\cases{
1-6x^2+6|x|^3,&\quad if  $|x|\leq\frac{1}{2}$,\cr
 2(1-|x|)^3, & \quad if\vspace*{2pt} $\frac{1}{2}\leq|x|\leq1$,\cr
  0, & \quad if $|x|>1$.
}
\end{equation}
Our analysis does not cover nontruncated kernels such as the quadratic
spectral kernel. But truncated kernels have the advantage of being
computationally more efficient.

Let $\{b_n,n\geq1\}$ be a nonincreasing sequence of positive numbers
such that
%
\begin{equation}\label{assumpSeq1}
b_n^{-1}=O(n^{1/2})\quad\mbox{and}\quad
|b_n-b_{n-1}|=O(b_nn^{-1})\qquad\mbox{as }n\to\infty.
\end{equation}

We consider the class of kernel estimator of the form
%
\begin{equation}
\Gamma^2_n(h)=\sum_{k=-n}^nw(k
b_n)\gamma_n(k)=\sum_{k=-b_n^{-1}+1}^{b_n^{-1}-1}w(k
b_n)\gamma_n(k).
\end{equation}

The following is the main $L^p$-convergence result.
\begin{theo}\label{thm1}
Assume \textup{A1}--\textup{A3}. Let $\beta\in
(0,1/2-\eta)$ and $h\in\mathcal{L}_{V^\beta}$, where $\eta$ is as in \textup{A2}.
Then
%
\begin{equation}\label{conclusionlem}
\Gamma^2_n(h)=\frac{1}{n}\sum_{k=1}^nG_{\theta_{k-1}}^2(X_{k-1},X_k)+
Q_n +D_n + \epsilon_n,\qquad n\geq1.
\end{equation}
The random process
$\{(Q_n,D_n,\epsilon_n),n\geq1\}$ is such that for any $p>1$ such
that $2p(\beta+\eta)\leq1$, there exists a finite constant $C$ such
that
\begin{eqnarray}\label{conclusionthm1}
\E(|Q_n|^p)&\leq& C\bigl(b_n + n^{-\alpha}b_n^{-1+\alpha} + n^{-1+
(1/2)\vee (1/p)} b_n^{-1/2}\bigr)^p,\nonumber\\ [-8pt]\\ [-8pt]
\E(|D_n|^p)&\leq& Cb_n^p \quad\mbox{and}\quad\E(|\epsilon_n|^p)\leq
C(n^{-1}b_n^{-1})^p.\nonumber
\end{eqnarray}
In particular, if $\lim_{n\to\infty}n^{-1+(1/2)\vee (1/p)}
b_n^{-1/2}=0,$ then
\[
\Gamma^2_n(h)-\frac{1}{n}\sum_{k=1}^nG_{\theta_{k-1}}^2 (X_{k-1}, X_k)
\]
converges to zero in $L^p$.\vspace*{-1pt}
\end{theo}
\begin{pf}
The proof is given in the supplementary article \cite{atchadesuppl}.\vspace*{-1pt}
\end{pf}

\begin{remark}
In Theorem \ref{thm1}, we can always take $p=1/(2(\beta+\eta))>1$. In
this case, the condition $\lim_{n\to\infty}n^{-1+(1/2)\vee
(1/p)} b_n^{-1/2}=0$ translates to $0.5\vee
(2(\beta+\eta))+0.5\delta<1$. Therefore, if $\beta+\eta$ is close to
$1/2$, we need to choose $\delta$ small. This remark implies that in
applying the above result, one should always try to find the smallest
possible $\beta$ such that $h\in\mathcal{L}_{V^\beta}$.

It can be easily checked that the choice of bandwidth $b_n\propto
n^{-\delta}$ with $\delta=\frac{2}{3}(1-0.5\vee
(2(\beta+\eta)))$ always satisfies Theorem \ref{thm1}. In fact,
we will see in Section \ref{remrate} that this choice of $b_n$ is
optimal in the $L^p$-norm, $p=(2(\beta+\eta))^{-1}$.

It is possible to investigate more carefully the rate of convergence of
$\Gamma_n^2(h)$ in Theorem \ref{thm1}. Indeed,\vspace*{1.5pt} consider the typical
case where $p=2$ is admissible and we have $\alpha=1$. If we choose
$b_n$ such that $b_n=o(n^{-1/3})$ and $n^{-1}=o(b_n),$ then the slowest
term in (\ref{conclusionthm1}) is $n^{-1+(1/2)\vee (1/p)}
b_n^{-1/2}=(n b_n)^{-1/2}$. By inspecting the proof of Theorem
\ref{thm1}, the only term whose $L^p$-norm enjoys such rate
$n^{-1+(1/2)\vee (1/p)} b_n^{-1/2}$ is
\[
Q_n^{(1)}=2n^{-1}\sum_{j=2}^n Z^{(1)}_{n,j}
G_{\theta_{j-1}}(X_{j-1},X_j),
\]
where
\[
Z_{n,j}^{(1)}=\sum_{\ell=1}^{j-1}w\bigl((j-\ell)b_n\bigr)G_{\theta_{\ell
-1}}(X_{\ell-1},X_\ell).
\]
Now, $\{(Q_n^{(1)},\F_n),n\geq2\}$ is a martingale array and we
conjecture that as $n\to\infty$,
\[
(n b_n)^{1/2}\Biggl(\Gamma^2_n(h)-\frac{1}{n}\sum_{k=1}^nG_{\theta
_{k-1}}^2(X_{k-1},X_k)\Biggr)
\stackrel{w}{\to}\mathcal{N}(0,\Lambda^2),\vadjust{\goodbreak}
\]
at least in the special
case where $\theta_n$ converges to a deterministic limit. But we do not
pursue this further since the issue of a central limit theorem for
$\Gamma_n^2(h)$ is less relevant for Monte Carlo simulation.
\end{remark}

When $\{X_n,n\geq0\}$ is a Markov chain, Theorem
\ref{thm1} improves on \cite{jonesetal09}, as it imposes weaker moment
conditions. Almost sure convergence is often more desirable in Monte
Carlo settings, but typically requires stronger assumptions. One can
impose either more restrictive growth conditions on $h$ (which
translates into stronger moment conditions, as in \cite{jonesetal09})
or one can impose stronger smoothness conditions on the function $w$.
We prove both types of results.

\begin{theo}\label{thm21}
Assume \textup{A1}--\textup{A3} with $\eta
<1/4$, where $\eta$ is as in \textup{A2}. Let $\beta\in(0,1/4-\eta)$
and $h\in\mathcal{L}_{V^\beta}$. Suppose that $b_n\propto n^{-\delta}$,
where $\delta\in(2(\beta+\eta),1/2)$. Then
\[
\lim_{n\to\infty}
\Biggl(\Gamma^2_n(h)-\frac{1}{n}\sum_{k=1}^nG_{\theta_{k-1}}^2(X_{k-1},X_k)\Biggr)=0
\]
almost surely.
\end{theo}
\begin{pf}
The proof is given in the supplementary article \cite{atchadesuppl}.
\end{pf}

We can remove the growth condition $h\in\mathcal{L}_{V^\beta}$,
$0<\beta<0.25-\eta,$ and the constraint on $b_n$ in Theorem \ref{thm21}
if we are willing to impose a stronger smoothness condition on $w$. To
do so, we replace A3 with A4.

\begin{longlist}
\item[A4:] The function $w$ is even [$w(-x)=w(x)$] and
$w(0)=1$. Moreover, the restriction $w\dvtx [0,1]\to\Real$ is $(r+1)$-times
continuously differentiable for some $r\geq2$.
\end{longlist}

\begin{theo}\label{thm2}
Assume \textup{A1}--\textup{A2} and \textup{A4}.
Let $\beta\in(0,1/2-\eta)$ and $h\in\mathcal{L}_{V^\beta}$, where $\eta$
is as in \textup{A2}. Let $p>1$ be such that $2p(\beta+\eta)\leq1$.
Suppose, in addition, that
\begin{eqnarray}\label{condSeq2}\sum_{n\geq
1}(n^{-1}b_n^{-1})^p&<&\infty,\qquad\sum_{n\geq1}(n^{-2}b_n^{-1})^{1\wedge
(p/2)}<\infty,\nonumber\\ [-8pt]\\ [-8pt]
\sum_{n\geq
1}n^{-2+(1/2)\vee (1/p)}b_n^{-1/2}&<&\infty\quad\mbox{and}\quad\sum_{n\geq1}
b_n^{(r-1)p}<\infty.\nonumber
\end{eqnarray}
The conclusion
of Theorem \ref{thm21} then holds.
\end{theo}
\begin{pf}
The proof is given in the supplementary article \cite{atchadesuppl}.
\end{pf}

\begin{remark}
Not all kernels used in practice will satisfy A4. For instance,
A4 holds for kernels in the family (\ref{bartlettsK}) but fails
to hold for the Parzen kernel (\ref{parzensK}).

In Theorem \ref{thm2}, we can again choose $b_n\propto n^{-\delta}$,
where $\delta=\frac{2}{3}(1-0.5\vee(2(\beta+\eta)))$. It is
easy to check that if A4 holds with
$r>1+2(\beta+\eta)\delta^{-1}$ and we take $p=(2(\beta+\eta))^{-1}$,
then this choice of $b_n$ satisfies (\ref{condSeq2}).
\end{remark}

In the next corollary, we consider the Markov chain case.

\begin{coro}\label{cor1}
Suppose that $\{X_n,n\geq0\}$ is a phi-irreducible, aperiodic Markov
chain with transition kernel $P$ and invariant distribution $\pi$.
Assume that $P$ satisfies \textup{A1}. Let $\beta\in(0,1/2)$ and
$h\in\mathcal{L}_{V^\beta}$. Then $\sigma^2(h):=\pi(h^2)+2\sum_{j\geq
1}\pi(h P^j h)$ is finite. Assume \textup{A3} and take
$b_n\propto n^{-\delta}$ with
$\delta=\frac{2}{3}(1-0.5\vee(2\beta))$. Then
\[
\lim_{n\to\infty} \Gamma^2_n(h)=\sigma^2(h)\qquad\mbox{in }L^{(2\beta)^{-1}}.
\]
Supposing, in addition, that
$\beta\in(0,1/4)$ and $\delta\in(2\beta,1/2),$ or that \textup{A4}
holds with $r>1+2\beta\delta^{-1}$, then the convergence holds almost
surely ($\p$) as well.
\end{coro}

\subsection{Application to the adaptive Metropolis algorithm}\label{sec41}
We shall now apply the above result to the adaptive Metropolis
algorithm described in Example \ref{exRWM}. We continue to use the
notation established in that example. We recall that $\Xset=\rset^d$,
$\Theta=\Theta_\mu(r_1)\times\Theta_+(r_2),$ where $\Theta_\mu
(r_1)$ is
the ball in $\Xset$ with center $0$ and radius $r_1>0$ and
$\Theta_+(r_2)$ is the set of all symmetric positive semidefinite
matrices $A$ with $|A|\leq r_2$. Define $\ell(x)=\log\pi(x)$. We assume
that:

\begin{longlist}
\item[B1:] $\pi$ is positive and continuously
differentiable,
\[
\lim_{|x|\to\infty}\biggl\langle\frac{x}{|x|},\nabla\ell(x)\biggr\rangle=-\infty
\]
and
\[
\lim_{|x|\to\infty}\biggl\langle\frac{x}{|x|},\frac{\nabla\ell
(x)}{|\nabla\ell(x)|}\biggr\rangle<0,
\]
where $\nabla\ell$ is the gradient of $\ell$.
\end{longlist}

B1 is known to imply A1 with
$V(x)=(\sup_{x\in\Xset}\pi^{\zeta}(x))\pi^{-\zeta}(x)$, for
any $\zeta\in(0,1)$ (\cite{andrieuetal06,jarnerethansen98}). We denote
by $\mu_\star$ and $\Sigma_\star$ the mean and covariance matrix of
$\pi$, respectively. We assume that
$(\mu_\star,\Sigma_\star)\in\Theta,$ which can always be achieved by
taking $r_1,r_2$ large enough.

By Lemma 12 of \cite{andrieuetal06}, for any $\beta\in(0,1]$,
%
\begin{equation}\label{diminishRWM}
D_\beta(\theta_n,\theta
_{n-1})\leq C|\Sigma_n-\Sigma_{n-1}|\leq\gamma_n V^\eta(X_n)
\end{equation}
for any $\eta>0$. Thus, A2 holds and $\eta$ can be taken to be
arbitrarily small. We can now summarize Proposition \ref{propclt} and
Theorems \ref{thm1}--\ref{thm2} for the random Metropolis algorithm. We
focus here on the choice of bandwidth $b_n\propto n^{-\delta}$, where
$\delta=\frac{2}{3}(1-0.5\vee(2\beta)),$ but similar
conclusions can be derived from the theorems for other bandwidths.

\begin{proposition}\label{propRWM}
Assume \textup{B1}, let
$V(x)=(\sup_{x\in\Xset}\pi^{\zeta}(x))\pi^{-\zeta}(x)$ for
$\zeta\in(0,1)$ and suppose\vadjust{\goodbreak} that $(\mu_\star,\Sigma_\star)\in
\Theta$.
Then $\theta_n=(\mu_n,\Sigma_n)$ converges in probability to
$\theta_\star=(\mu_\star,\Sigma_\star)$. Let $\beta\in(0,1/2)$ and
$h\in\mathcal{L}_{V^\beta}$.
\begin{longlist}
\item[1.]$n^{-1/2}\sum_{k=1}^n h(X_k)$ converges weakly to
$\mathcal{N}(\pi(h),\sigma^2_\star(h))$ as $n\to\infty$,
where $\sigma^2_\star(h)=\pi(h^2)+2\sum_{j\geq
1}\pi(hP_{\theta_\star}^jh)$ and
$\theta_\star=\Sigma_\star+\epsilon I_d$.
\item[2.] Suppose that \textup{A3}
holds and we choose $b_n\propto n^{-\delta}$,
$\delta=\frac{2}{3}(1-0.5\vee(2\beta))$. Then
$\Gamma_n^2(h)$ converges to $\sigma^2_\star(h)$ in $L^p$ for
$p=(2\beta)^{-1}$. If we \vspace*{1pt}additionally suppose that $\beta\in(0,1/4)$
and $\delta\in(2\beta,1/2),$ or that \textup{A4} holds with
$r>1+2\beta\delta^{-1}$, then the convergence of $\Gamma_n^2(h)$ holds
almost surely ($\p$) as well.
\end{longlist}
\end{proposition}

\subsection{Choosing the bandwidth $b_n$}\label{remrate}
Consider Theorem \ref{thm1}. Suppose that $\alpha\geq2/3$ and that we
take $b_n\propto n^{-\delta}$ for some $\delta\in(0,1/2]$. Then
$n^{-\alpha}b_n^{-1+\alpha}=O(n^{-1/2})$. Similarly,
$n^{-1}b_n^{-1}=O(n^{-1/2})$. Thus, the $L^p$-rate of convergence of
$\Gamma^2_n(h)$ is driven by $b_n$ and $n^{-1+(1/2)\vee
(1/p)} b_n^{-1/2}$, and we deduce from equating these two\vspace*{1.5pt} terms
that the optimal choice of $b_n$ is given by $b_n \propto n^{-\delta}$
for $\delta=\frac{2}{3}(1-\frac{1}{2}\vee\frac{1}{p})$.
Equation~(\ref{conclusionthm1}) then gives that\vspace{-3pt}
\[
\E^{1/p}\Biggl(\Biggl|\Gamma^2_n(h)-\frac{1}{n}\sum_{k=1}^nG_{\theta
_{k-1}}^2(X_{k-1},X_k)\Biggr|^p\Biggr)\leq Cn^{-\delta}.
\]
In particular, if $4(\beta+\eta)\leq1$ (and $\alpha\geq2/3$), we can
take $p=2$ and then $\delta=1/3,$ which leads to\vspace{-3pt}
\[
\E^{1/2}\Biggl(\Biggl|\Gamma^2_n(h)-\frac{1}{n}\sum_{k=1}^nG_{\theta
_{k-1}}^2(X_{k-1},X_k)\Biggr|^2\Biggr)\leq Cn^{-1/3}.
\]
The same $L^2$-rate of convergence was also derived in
\cite{jonesetal09}.

Even with $b_n=\frac{1}{cn^{1/3}}$, the estimator is still very
sensitive \vspace*{1.5pt}to the choice of~$c$. Choosing $c$ is a difficult issue where
more research is needed. Here, we follow a data-driven approach adapted
from \cite{andrew91} and \cite{neweyetwest94}. In this approach, we
take $b_n=\frac{1}{cn^{1/3}}$, where\vspace{-5pt}
\[
c=c_0\biggl\{\frac{2\sum_{\ell=1}^m\ell\hat\rho_\ell}{1+2\sum_{\ell
=1}^m\hat\rho_\ell}\biggr\}^{1/3}
\]
for some constants $c_0$ and $m$, where $\hat\rho_\ell$ is the $\ell$th
order sample autocorrelation of $\{h(X_n),n\geq0\}$.
\cite{neweyetwest94} suggests choosing $m=n^{2/9}$. Our simulation
results show that small values of $c_0$ yield small variances but high
biases, and inversely for large values of $c_0$. The value $c_0$ also
depends on how fast the autocorrelation of the process decays.
\cite{neweyetwest94} derives some theoretical results on the
consistency of this procedure in the stationary case. Whether these
results hold in the present nonstationary case is an open question.

\subsection{Discussion}\label{finalrem}
The above results raise a number of issues. On one hand, we note from
Theorems \ref{thm1}--\ref{thm2} that the kernel estimator
$\Gamma_n^2(h)$ does not converge to the asymptotic variance
$\sigma^2(h)$, but rather to the asymptotic average squared variation
$\Gamma^2(h)$. On the other hand, Proposition \ref{propclt} shows that
although the asymptotic variance $\sigma^2(h)$ controls the
fluctuations of $n^{-1/2}\sum_{k=1}^k h(X_k)$ as $n\to\infty$, the
limiting distribution of $n^{-1/2}\sum_{k=1}^k h(X_k)$ is not the
Gaussian $\mathcal{N}(0,\sigma^2(h))$, but instead a mixture of
Gaussian distribution of the form $\sqrt{\Gamma^2(h)}Z$. With these
conditions, how can one undertake a valid error assessment from
adaptive MCMC samplers?

If the adaptation parameter $\theta_n$ converges to a deterministic
limit $\theta_\star$, then one gets a situation similar to that of
Markov chains. This is the ideal case. Indeed, in such cases,
$\Gamma^2(h)\equiv\sigma^2(h)$, $n^{-1/2}\sum_{k=1}^nh(X_k)$ converges
weakly to a random variable $\mathcal{N}(0,\sigma^2(h))$ and
the kernel estimator $\Gamma_n^2(h)$ converges to the asymptotic
variance $\sigma^2(h)$, where
\[
\sigma^2(h)=\int_{\Xset\times\Xset}\pi(dx)P_{\theta_\star
}(x,dy)G_{\theta_\star}^2(x,y)=\pi(h^2)+2\sum_{j\geq1}\pi
(hP_{\theta_\star}^jh).
\]
This case includes the adaptive Metropolis algorithm of
\cite{haarioetal00}, as discussed in Section \ref{sec41}.

However, in some other cases (see, e.g.,
\cite{andrieuetal06,atchadeetfort08}), what one can actually prove is
that $\theta_n\to\theta_\star$, where $\theta_\star$ is a discrete
random variable with values in a subset
$\{\tau_1,\tau_2,\ldots,\tau_N\},$ say, of $\Theta$. This is typically
the case when the adaptation is driven by a stochastic approximation
$\theta_{n+1}=\theta_n+\gamma_n H(\theta_n,X_{n+1}),$ where the mean
field equation $h(\theta)\eqdef\int_\Xset H(\theta,x)\pi(dx)=0$ has
multiple \vspace*{1pt}solutions.

In these cases, $\Gamma^2_n(h)$ clearly provides a poor estimate for
$\sigma^2(h)$, even though it is not hard to see that
\[
\lim_{n\to\infty}\E(\Gamma_n^2(h))=\PE(\Gamma^2(h))=\sigma^2(h).
\]
Furthermore, a confidence interval for $\pi(h)$ becomes difficult to
build. Indeed, the asymptotic distribution $n^{-1/2}\sum_{k=1}^n
h(X_k)$ is a mixture
\[
\sum_{k\geq1}p_k\mathcal{N}(0,\sigma_k),
\]
where $p_k\eqdef\p(\theta_\star=\tau_k)$ and
$\sigma^2_k(h)=\pi(h^2)+2\sum_{j\geq1}\pi(hP_{\tau_k}^jh)$.
As\vspace*{1pt} a consequence, a valid confidence interval for $\pi(h)$ requires
the knowledge\vspace*{1pt} of the mixing distribution $p_k$ and the asymptotic
variances $\sigma_k^2(h)$, which is much more than one can obtain from
$\Gamma_n^2(h)$. It is possible to improve on the estimation of
$\sigma^2(h)$ by running multiple chains, but this takes away some of
the advantages of the adaptive MCMC framework.

In view of this discussion, when Monte Carlo error assessment is
important, it seems that the framework of adaptive MCMC is most useful
when the adaptation mechanism is such that there exists a unique,
well-defined, optimal kernel $P_{\theta\star}$ that the algorithm
converges to. This is the case, for example, with the popular adaptive
RWM of \cite{haarioetal00} discussed above and its extension to the
MALA (Metropolis adjusted Langevin algorithm; see, e.g.,~\cite{atchade05}).

\section{Examples}\label{sec:examples}
\subsection{The $\operatorname{GARCH}(1,1)$ model}

To illustrate the above results in the Markov chain case, we consider
the linear $\operatorname{GARCH}(1,1)$ model defined as follows: $h_0\in(0,\infty)$,
$u_0\sim\mathcal{N}(0,h_0)$ and, for $n\geq1,$
\begin{eqnarray*}
u_n&=&h_n^{1/2}\epsilon_n,\\
h_n&=&\omega+\beta h_{n-1} + \alpha u_{n-1}^2,
\end{eqnarray*}
where
$\{\epsilon_n,n\geq0\}$ is i.i.d. $\mathcal{N}(0,1)$ and $\omega>0$,
$\alpha\geq0$, $\beta\geq0$. We assume that $\alpha,\beta$ satisfy
the following.

\begin{longlist}
\item[E1:] There exists $\nu>0$ such that
%
\begin{equation}\label{E1}
\E[(\beta+\alpha
Z^2)^\nu]<1,\qquad Z\sim\mathcal{N}(0,1).
\end{equation}
\end{longlist}

It is shown by \cite{meitzetsaikkonen08}, Theorem 2, that under
(\ref{E1}), the joint process $\{(u_n,h_n),\break n\geq0\}$ is a
phi-irreducible aperiodic Markov chain that admits an invariant
distribution and is geometrically ergodic with a drift function
$V(u,h)=1+h^\nu+|u|^{2\nu}$. Therefore, A1 holds and we can
apply Corollary \ref{cor1}. We write $\E_\pi$ to denote expectation
taken under the stationary measure. We are interested in the asymptotic
variance of the functions $h(u)=u^2$. We can calculate the exact value.
Define $\rho_n\eqdef\operatorname{Corr}_\pi(u_0^2,u_n^2)$. As observed by
\cite{bollerslev86} in introducing the $\operatorname{GARCH}$ models,
if (\ref{E1}) hold with some $\nu\geq2$, then
\[
\rho_1=\frac{\alpha(1-\alpha\beta-\beta^2)}{1-2\alpha\beta
-\beta^2},\qquad\rho_n=\rho_1(\alpha+\beta)^{n-1},\qquad n\geq2.
\]
Also,
\[
\operatorname{Var}_\pi(u_0^2)=\frac{3\omega^2(1+\alpha+\beta)}{(1-\alpha
-\beta)(1-\beta^2-2\alpha\beta-3\alpha^2)}-\biggl(\frac{\omega
}{1-\alpha-\beta}\biggr)^2
\]
and we obtain
\[
\sigma^2(h)=\operatorname{Var}_\pi(u_0^2)\biggl(1+2\frac{\rho_1}{1-\alpha-\beta}\biggr).
\]

\begin{figure}

\includegraphics{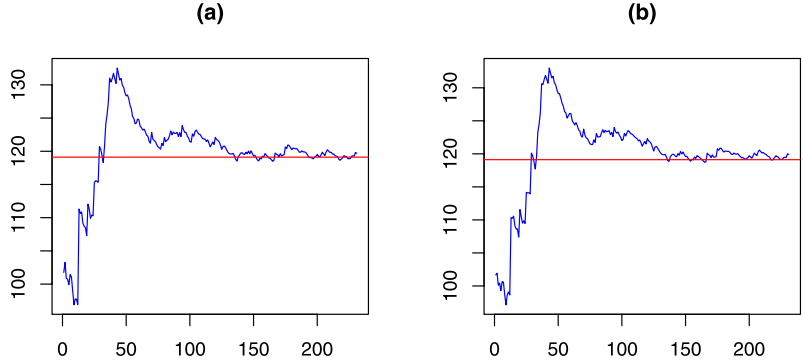}

\caption{Asymptotic variance estimation for
$\operatorname{GARCH}(1,1)$ with $\omega=1$, $\alpha=0.1$, $\beta=0.7$ based on
250,000 iterations. \textup{(a)} is Bartlett kernel, \textup{(b)} is Parzen
kernel.}\label{fig1}
\end{figure}

For the simulations, we set $\omega=1$, $\alpha=0.1$, $\beta=0.7,$
which gives $\sigma^2(h)=119.1$. For these values, (\ref{E1}) holds
with at least $\nu=4$. We tested the Bartlett and the Parzen kernels
for which A3 holds. We choose the bandwidth following the
approach outlined in Remark \ref{remrate} with $c_0=1.5$. We run the
$\operatorname{GARCH}(1,1)$ Markov chain for 250,000 iterations and discard the first
10,000 iterations as burn-in. We compute $\Gamma^2_n(h)$ at every 1000
along the sample path. The results are plotted in Figure \ref{fig1}.

\subsection{Logistic regression}
We also illustrate the results with MCMC and adaptive MCMC. We consider
the logistic regression model
\[
y_i\sim\mathcal{B}(p_\beta(x_i)),\qquad i=1,\ldots,n,
\]
where $y_i\in\{0,1\}$ and $p_\beta(x)=e^{x\beta}(1+e^{x\beta})^{-1}$
for a parameter $\beta\in\Real^d$ and a covariate vector
$x^T\in\Real^d$, where $x^T$ denotes the transpose of $x$.
$\mathcal{B}(p)$ is the Bernoulli distribution with parameter $p$. The
log-likelihood is
\[
\ell(\beta\vert X)=\sum_{i=1}^ny_ix_i\beta
-\log(1+e^{x_i\beta}).
\]
We assume a Gaussian prior
distribution $\pi(\beta)\propto e^{-1/(2s^2)|\beta|^2}$ for some
constant $s>0$ leading to a posterior distribution
\[
\pi(\beta|X)\propto e^{\ell(\beta\vert
X)}e^{-1/(2s^2)|\beta|^2}.
\]
The RWM algorithm described in Example \ref{exRWM} is a possible choice to
sample from the posterior distribution. We compare a plain RWM with
proposal density $\mathcal{N}(0,e^cI_d)$ with $c=-2.3$ and the adaptive
RWM described in Algorithm \ref{arwm1} using the family
$\{P_\theta,\theta\in\Theta\},$ where $\Theta=\Theta_\mu(r_1)
\times
\Theta_+(r_2),$ as defined in Example~\ref{exRWM}. It is easy to check that
B1 holds. Indeed, we have
\[
\langle\beta,\nabla\log\pi(\beta)\rangle=-\frac{|\beta|^2}{s^2}+\sum
_{i=1}^n\bigl(y_i-p_\beta(x_i)\bigr)\langle\beta,x_i^T\rangle
\]
and $|\sum_{i=1}^n(y_i-p_\beta(x_i))\langle\beta,x_i^T\rangle|\leq
|\beta|\sum_{i=1}^n|x_i|$. We deduce that
\[
\biggl\langle\frac{\beta}{|\beta|,\nabla\log\pi(\beta)}\biggr\rangle\leq-\frac
{|\beta|}{s^2}+\sum_{i=1}^n|x_i|\to-\infty\qquad\mbox{as }|\beta|\to
\infty.
\]

Similarly,
\[
\biggl\langle\frac{\beta}{|\beta|},\frac{\nabla\log\pi(\beta)}{|\nabla
\log\pi(\beta)|}\biggr\rangle
\leq-\frac{1}{s^2}\frac{|\beta|}{|\nabla\log\pi(\beta)|}+\frac
{\sum_{i=1}^n|x_i|}{|\nabla\log\pi(\beta)|}
\to-1\qquad\mbox{as }|\beta|\to\infty
\]
since
$|\nabla\log\pi(\beta)|\sim s^{-2}|\beta|$ as $|\beta|\to\infty$.
Therefore, B1 holds. If we choose $r_1,r_2$ large enough so that
$(\mu_\star,\Sigma_\star)\in\Theta,$ then Proposition \ref{propRWM}
holds and applies to any measurable function $h$ such that
$|h(\beta)|\leq c \pi^{-t}(\beta|X)$ for some $t\in[0,1/2)$.
\begin{figure}

\includegraphics{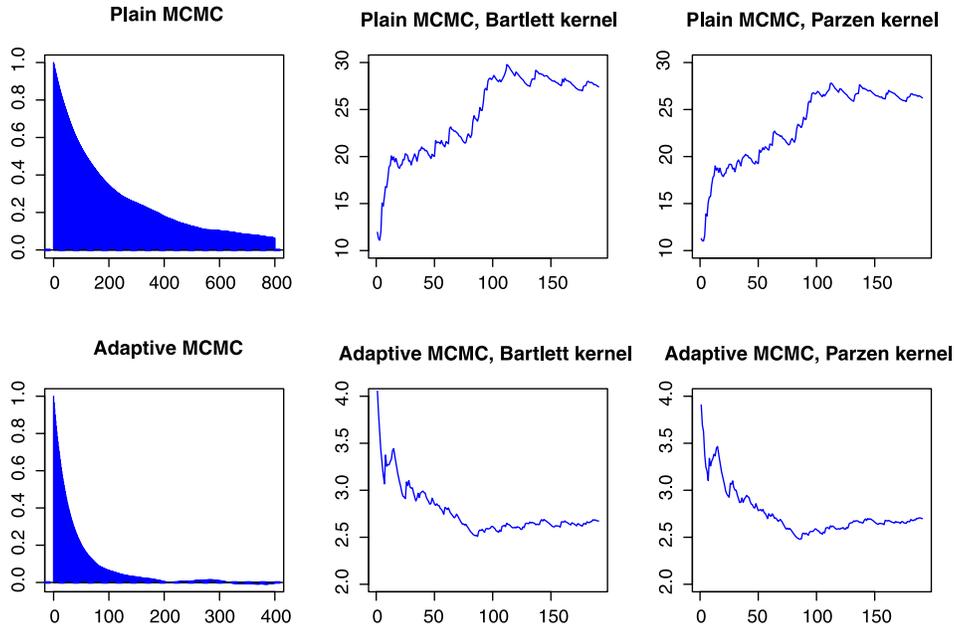}

\caption{Asymptotic variance estimation for logistic
regression modeling of the heart data set. Outputs of the coefficient
$\beta^{(2)}$ are reported, based on 250,000 iterations.}\label{fig2}
\end{figure}

As a simulation example, we test the model with the \textit{Heart data
set} which has $n=217$ cases and $d=14$ covariates. The dependent
variable is the presence or absence of a heart disease and the
explanatory  variables are relevant covariates. More details can be found in
\cite{michieetal94}. We use Parzen and Bartlett kernels with $c_0=20$
for the Markov chain and $c_0=5$ for the adaptive chain. We run both
chains for 250,000 iterations and discard the first 50,000 iterations
as burn-in. The results are plotted in Figure \ref{fig2} for the coefficient
$\beta_2$. We also report in Table \ref{table1} below the resulting confidence for
the first four coefficients $(\beta_1,\ldots,\beta_4)$.

\section{Proofs}\label{sec:proofs}

This section contains the proofs of the statements from Sections \ref{sec:adaptchain}--\ref{sec:clt}. The remaining proofs are available
in the supplementary paper \cite{atchadesuppl}. Throughout this
section, we shall use $C$ to denote a generic constant whose actual
value might change from one appearance to the next. On multiple
occasions, we make use of the Kronecker lemma and the Toeplitz lemma.
We refer the reader to \cite{halletheyde80}, Section~2.6, for a
statement and proof of these lemmata.

We shall routinely use the following martingale inequality. Let
$\{D_{i},\F_{i},i\geq1\}$ be a martingale difference sequence. For
any $p>1$,
%
\begin{equation}\label{boundM}
\E\Biggl(\Biggl|\sum_{i=1}^n D_i\Biggr|^p\Biggr)\leq
C\Biggl\{\sum_{i=1}^n\PE^{1\wedge(2/p)}(|D_{i}|^p)\Biggr\}^{1\vee
(p/2)},\vspace*{1pt}
\end{equation}
where $C$ can be taken as
$C=(18pq^{1/2})^p$, $p^{-1}+q^{-1}=1$.

\begin{table}[t]
\tablewidth=235pt
\tabcolsep=0pt
\caption{Confidence interval for the first four parameters of the model
for the heart data set}\label{table1}
\begin{tabular*}{235pt}{@{\extracolsep{\fill}}lcc@{}}
\hline
\textbf{Parameters}&\textbf{Plain MCMC}&\textbf{Adaptive RWM}\\
\hline
$\beta_1$ & $[-0.271,-0.239]$ & $[-0.272,-0.257]$\\
$\beta_2$ & $[-0.203,-0.158]$ & $[-0.182,-0.170]$\\
$\beta_3$ & $[0.744,0.785]$& $[0.776,0.793]$\\
$\beta_4$ & $[0.727,0.756]$ & $[0.736,0.750]$\\
\hline
\end{tabular*}\vspace*{-4pt}%
\end{table}

We also notice that for any $q\in[1,\beta^{-1}],$ Lemma
\ref{lemuseful}(a)--(b) implies that
%
\begin{equation}\label{momentG}\sup_{k\geq1} \E(|G_{\theta
_{k-1}}(X_{k-1},X_k)|^{q})<\infty.\vspace*{-4pt}
\end{equation}

\subsection{\texorpdfstring{Proof of Proposition \protect\ref{wlln}}{Proof of Proposition 2.1}}
\label{proof:propslln}

Let $S_n\eqdef\sum_{k=1}^n h_{\theta_{k-1}}(X_k)$. For
$\theta\in\Theta$, we define $\tilde g_\theta(x)=\sum_{j\geq0}
P_\theta^j h_\theta(x)$. When $h_\theta$ does not depend on $\theta$,
we obtain $\tilde g_\theta=g_\theta$, as defined in Section \ref{sec:clt}.
Similarly, we define $\tilde G_\theta(x,y)=\tilde g_\theta
(y)=P_\theta
\tilde g_\theta(x)$. Using the Poisson equation $\tilde g_\theta-
P_\theta\tilde g_\theta=h_\theta$, we rewrite $S_n$ as $S_n=M_n +R_n$,
where
\[
M_n\eqdef\sum_{k=1}^n \tilde G_{\theta_{k-1}}(X_{k-1},X_k)
\]
and
\[
R_n\eqdef P_{\theta_0}\tilde g_{\theta_0}(X_0)-P_{\theta_n}\tilde
g_{\theta_n}(X_n)+\sum_{k=1}^n\bigl(\tilde g_{\theta_{k}}(X_k)-\tilde
g_{\theta_{k-1}}(X_k)\bigr).
\]

Using Lemma \ref{lemuseful} and A1, we easily see that
\[
|R_n|\leq C\Biggl(V^\beta(X_0)+V^\beta(X_n)+\sum_{k=1}^n\bigl(D_\beta(\theta
_k,\theta_{k-1})+|h_{\theta_k}-h_{\theta_{k-1}}|_{V^\beta}\bigr)V^{\beta}(X_k)\Biggr).
\]
For $p>1$ such that $\beta p\leq1$, $\sum_{k\geq
1}n^{-p}\E((V^\beta(X_0)+V^\beta(X_n))^p)<\infty$. This is~a~consequence
of Lemma \ref{lemuseful}(a) and the Minkowski inequality.
Thus, $n^{-1}\times(V^\beta(X_0)+ V^\beta(X_n))$ converges almost surely to
zero. By (\ref{diminishwlln}) and\vspace*{1pt} the Kronecker lemma, the term
$n^{-1}\sum_{k=1}^n(D_\beta(\theta_k,\theta_{k-1})+|h_{\theta
_k}-h_{\theta_{k-1}}|_{V^\beta})V^{\beta}(X_k)$
converges almost surely to zero. We conclude that $n^{-1}R_n$ converges
almost surely to zero.\vadjust{\goodbreak}

$\{(M_n,\F_n),n\geq1\}$ is a martingale. Again, let $p>1$ be\vspace*{1pt} such
that $\beta p\leq1$.
Equation (\ref{boundM}) and Lemma \ref{lemuseful}(a)
together imply that $\PE(|M_n|^p)=O(n^{1\vee
(p/2)})$, which, combined with Proposition A.1 of
\cite{atchadesuppl}, implies that $n^{-1}M_n$ converges almost surely
to zero.

\subsection{\texorpdfstring{Proof of Proposition \protect\ref{propclt}}{Proof of Proposition 3.1}}
\label{proof:propclt}
This is a continuation of the previous\vspace*{1.5pt} proof. In the
present case, $h_\theta\equiv h$, so we write $g_\theta$ and
$G_\theta$ instead\vspace*{-1pt} of $\tilde g_\theta$ and $\tilde G_\theta$,
respectively. Again, let $S_n\eqdef\sum_{k=1}^n h(X_k)$. We have
$S_n=M_n +R_n$, where $M_n\eqdef\sum_{k=1}^n G_{\theta
_{k-1}}(X_{k-1},X_k)$ and
\[
|R_n|\leq C\Biggl(V^\beta(X_0)+V^\beta(X_n)+\sum_{k=1}^nD_\beta(\theta
_k,\theta_{k-1})V^{\beta}(X_k)\Biggr).
\]
$\E(V^\beta(X_0)+V^\beta(X_n))$ is bounded in $n$, thus
$n^{-1/2}(V^\beta(X_0)+V^\beta(X_n))$ converges in probability to zero.
By (\ref{diminishclt}) and the Kronecker lemma, the term
$n^{-1/2}\sum_{k=1}^nD_\beta(\theta_k,\theta_{k-1})V^{\beta}(X_k)$
converges \vspace*{1.5pt}almost surely to zero.\vspace*{1pt} We conclude that $n^{-1/2}R_n$
converges in probability to zero.

$\{(M_n,\F_n),n\geq1\}$ is a martingale. Since $\beta<1/2$,
(\ref{momentG}) implies that $\{(M_n,\break\F_n), n\geq1\}$ is a square
integrable martingale and also that we have
\begin{eqnarray}\label{propeq1}
\sup_{n\geq1}\E\Bigl(\max_{1\leq k\leq n}
n^{-1}G^2_{\theta_{k-1}}(X_{k-1},X_k)\Bigr)&<&\infty\quad\mbox{and}\nonumber\\
[-8pt]\\ [-8pt]
\lim_{n\to\infty}\max_{1\leq k\leq n}
n^{-1/2}G_{\theta_{k-1}}(X_{k-1},X_k)&=&0 \qquad\mbox{(in
probability)}.\nonumber
\end{eqnarray}
Equations\vspace*{1pt} (\ref{asympvar}) and (\ref{propeq1})
imply, by Theorem 3.2 of \cite{halletheyde80}, that $n^{-1/2}M_n$
converges weakly \vadjust{\goodbreak}to a random variable $\sqrt{\Gamma^2(h)}Z$, where
$Z\sim\mathcal{N}(0,1),$ and is independent of $\Gamma^2(h)$.

\subsection{\texorpdfstring{Proof of Proposition \protect\ref{propconvtheta}}{Proof of Proposition 3.3}}
\label{proof:convtheta}
We have
\begin{eqnarray*}
&&\frac{1}{n}\sum_{k=1}^nG^2_{\theta_{k-1}}(X_{k-1},X_k)\\
&&\qquad=\frac{1}{n}\sum_{k=1}^n\bigl(G^2_{\theta_{k-1}}(X_{k-1},X_k)-P_{\theta
_{k-1}}G^2_{\theta_{k-1}}(X_{k-1})\bigr)\\
&&\qquad\quad{}+ \frac{1}{n}\sum_{k=1}^n\biggl(P_{\theta_{k-1}}G^2_{\theta
_{k-1}}(X_{k-1})-\int_\Xset\pi(dx) P_{\theta_{k-1}}G^2_{\theta
_{k-1}}(x)\biggr)\\
&&\qquad\quad{}+\frac{1}{n}\sum_{k=1}^n\int_\Xset\pi(dx) \bigl(P_{\theta
_{k-1}}G^2_{\theta_{k-1}}(x)-P_{\theta_{\star}}G^2_{\theta_{\star
}}(x)\bigr)+\int_\Xset P_{\theta_\star}G^2_{\theta_\star}(x)\pi(dx)\\
&&\qquad=T_n^{(1)}+T_n^{(2)}+T_n^{(3)}+\int_\Xset
P_{\theta_\star}G^2_{\theta_\star}(x)\pi(dx),
\end{eqnarray*}
say. The
term $T_n^{(1)}$ is an $\F_n$-martingale. Indeed,
$\PE(G^2_{\theta_{k-1}}(X_{k-1},X_k)\vert
\F_{k-1})=P_{\theta_{k-1}}G^2_{\theta_{k-1}}(X_{k-1})$,
$\PP$-a.s. Furthermore, by (\ref{momentG}), the martingale differences
$G^2_{\theta_{k-1}}(X_{k-1},X_k)-P_{\theta_{k-1}}G^2_{\theta_{k-1}}(X_{k-1})$
are $L^p$-bounded for some $p>1$. By~\cite{halletheyde80}, Theorem 2.22,
we conclude that $T_n^{(1)}$ converges in $L^1$ to zero.

The term $T_n^{(2)}$ converges in probability to zero as a consequence
of the law of large numbers (Proposition \ref{wlln}). Using the
definition of $D_\beta$ and Lemma \ref{lemuseful}(a)--(b), we can find
a constant $C$ such that
\begin{eqnarray*}
&&\biggl|\int_\Xset\pi(dx) \bigl(P_{\theta_{n}}G^2_{\theta_{n}}(x)-P_{\theta
_{\star}}G^2_{\theta_{\star}}(x)\bigr)\biggr|\\
&&\qquad\leq C \bigl(D_{\beta}(\theta_n,\theta_\star) +
D_{2\beta}(\theta_n,\theta_\star)\bigr)\int_\Xset
V^{2\beta}(x)\pi(dx),
\end{eqnarray*}
almost surely. It follows that
$T_n^{(3)}$ also converges in $\PP$-probability to zero.

\subsection{\texorpdfstring{Proof of Proposition \protect\ref{propvar}}{Proof of Proposition 3.4}}
\label{proof:propvar}
From the proof of Proposition \ref{wlln} above, we have seen that
$S_n=M_n+R_n,$ and it is easy to check that
$\E(|R_n|^2)=O(n^{2(1-\alpha)})$ and,
by (\ref{momentG}), $\E(|M_n|^2)=O(n)$.
Therefore,
\begin{eqnarray*}
&&|\operatorname{Var}(n^{-1/2}S_n)-n^{-1}\PE(M_n^2)|\\
&&\qquad=|2n^{-1}\E(M_n R_n) +n^{-1}\PE(R_n^2) -n^{-1}(\E(R_n))^2|\\
&&\qquad=O(n^{1/2-\alpha})\to0\qquad\mbox{as }n\to\infty
\end{eqnarray*}
since $\alpha>1/2$. Now,
\[
n^{-1}\E(M_n^2)=\E\Biggl(n^{-1}\sum_{k=1}^n G_{\theta_{k-1}}^2(X_{k-1},X_k)\Biggr).
\]
Again, from (\ref{momentG}), the sequence $n^{-1}\sum_{k=1}^n
G_{\theta_{k-1}}^2(X_{k-1},X_k)$ is uniformly integrable which,
combined with (\ref{asympvar}) and Lebesgue's dominated convergence
theorem, implies that $n^{-1}\E(M_n^2)$ converges to $\E(\Gamma^2(h))$.

\begin{appendix}

\section{Some useful consequences of A1}
\begin{lemma}\label{lemuseful}
Assume that $\{P_\theta,\theta\in\Theta\}$ satisfies \textup{A1}.
\begin{longlist}
\item[(a)] There exists a finite constant $C$ such that
%
\begin{equation}\label{lemusefuleq1}
\sup_{n\geq0}\PE(V(X_n))\leq C.
\end{equation}
\item[(b)]
Let $\beta\in(0,1]$ and
$\{h_\theta\in\mathcal{L}_{V^\beta},\theta\in\Theta\}$ be such that
$\pi(h_\theta)=0$, $\sup_{\theta\in\Theta}\break|h_\theta|_{V^\beta
}<\infty$.
The function $\tilde g_\theta\eqdef\sum_{j\geq0}P^j_\theta
h_\theta(x)$ is then well defined, $|\tilde
g_\theta|_{V^\beta}\leq C|h_\theta|_{V^\beta}$, where the
constant $C$ does not depend on
$\{h_\theta\in\mathcal{L}_{V^\beta},\theta\in\Theta\}$. Moreover, we can take
$C$ such that for any $\theta,\theta'\in\Theta$,
%
\begin{equation}\label{lemusefuleq2}
|\tilde g_\theta-\tilde g_{\theta
'}|_{V^\beta}\leq C\sup_{\theta\in\Theta}|h_\theta|_{V^\beta}
\bigl(D_\beta(\theta,\theta')+|h_\theta-h_{\theta'}|_{V^\beta}\bigr).
\end{equation}
\item[(c)] Assume \textup{A2}. Let $\beta\in(0,1-\eta)$ and
$h\in\mathcal{L}_{V^\beta}$ be such that
$\pi(h)=0$. Define $S_n(j)=\sum_{\ell=j+1}^{j+n} h(X_\ell)$. Let
$p\in(1,(\beta+\eta)^{-1})$. There then exists a~fini\-te constant $C$
that does not depend on $n,j,\theta$ or $h$
such that
\[
\PE(|S_n(j)|^p)\leq C |h|_{V^\beta}n^{1\vee(p/2)}.
\]
\end{longlist}
\end{lemma}
\begin{pf}
Parts (a) and (b) are standard results (see, e.g.,
\cite{andrieuetal06}). To prove~(c), we use the Poisson equation
(\ref{poisson}) to write
\begin{eqnarray*}
S_n(j)&=&\sum_{\ell=j+1}^{j+n}G_{\theta_{\ell
-1}}(X_{\ell-1},X_\ell)+P_{\theta_j}g_{\theta_j}(X_j)
-P_{\theta_{j+n}}g_{\theta_{j+n}}(X_{j+n})\\
&&{}+\sum_{\ell=j+1}^{j+n}\bigl(g_{\theta_{\ell-1}}(X_\ell)-g_{\theta_\ell}(X_l)\bigr).
\end{eqnarray*}
By A1 and part (a), we have
\[
\sup_{n\geq1}\sup_{j\geq
0}\E[|P_{\theta_j}g_{\theta_j}(X_j)
-P_{\theta_{j+n}}g_{\theta_{j+n}}(X_{j+n})|^p]\leq C
|h|_{V^\beta}.
\]
By Burkholder's inequality and some standard
inequalities,
\begin{eqnarray*}
\E\Biggl[\Biggl|\sum_{\ell=j+1}^{j+n}G_{\theta_{\ell-1}}(X_{\ell-1},X_\ell
)\Biggr|^p\Biggr]&\leq& C\Biggl\{\sum_{\ell=j+1}^{j+n}\E^{1\wedge(2/p)}(|G_{\theta
_{\ell-1}}(X_{\ell-1},X_\ell)|^p)\Biggr\}^{1\vee(p/2)}\\
&\leq& C|h|_{V^\beta}n^{1\vee(p/2)}.
\end{eqnarray*}
Part (b) and
A2 together give
\begin{eqnarray*}
&&\E\Biggl[\Biggl|\sum_{\ell=j+1}^{j+n}g_{\theta_{\ell-1}}(X_\ell)-g_{\theta
_\ell}(X_l)\Biggr|^p\Biggr]\\
&&\qquad\leq C|h|_{V^\beta}\E\Biggl[\Biggl(\sum_{\ell=j+1}^{j+n}D_\beta(\theta_{\ell-1},\theta_\ell)
V^{\beta}(X_\ell)\Biggr)^p\Biggr]\\
&&\qquad\leq C|h|_{V^\beta}\E\Biggl[\Biggl(\sum_{\ell=j+1}^{j+n}\gamma_{k+\ell}V^{\beta
+\eta}(X_\ell)\Biggr)^p\Biggr]\leq
C|h|_{V^\beta}\Biggl(\sum_{\ell=j+1}^{j+n}\gamma_{k+\ell}\Biggr)^p
\end{eqnarray*}
and, since $\gamma_n=O(n^{-1/2}),$ we are done.
\end{pf}
\end{appendix}
\section*{Acknowledgments}
The author is grateful to Galin Jones for helpful discussions,
and to Prosper Dovonon for pointing out some of the references in
the econometrics literature and for helpful comments on an earlier
version of this paper.

\begin{supplement}[id=suppA]
\stitle{Supplement to ``Kernel estimators of asymptotic variance for
adaptive Markov chain Monte Carlo''}
\slink[doi]{10.1214/10-AOS828SUPP}
\sdatatype{.pdf}
\sfilename{Suppl.pdf}
\sdescription{The proofs of Theorems \ref{thm1}--\ref{thm2} require
some technical and lengthy arguments that we develop in this supplement.}
\end{supplement}

\printaddresses

\end{document}